\documentclass[a4paper]{amsart}

\usepackage{graphicx}
\usepackage{subfigure}
\thanks{The first author was partially supported by The European
Contract Human Potential Programme, Research Training Network
HPRN-CT-2000-00101} 
\title{The Ribaucour transformation in Lie sphere geometry}
\author{F.E. Burstall}
\address{Department of Mathematical Sciences\\
University of Bath\\Bath BA2 7AY\\UK}
\email{feb@maths.bath.ac.uk}
\author{U. Hertrich-Jeromin}
\address{Department of Mathematical Sciences\\
University of Bath\\Bath BA2 7AY\\UK}
\email{masuh@maths.bath.ac.uk}

\newtheorem{thm}{Theorem}[section]
\newtheorem{lemma}[thm]{Lemma}
\newtheorem{prop}[thm]{Proposition}
\newtheorem*{BPT}{Bianchi Permutability Theorem}
\newtheorem*{CT}{Cube Theorem}

\theoremstyle{definition}
\newtheorem*{defn}{Definition}
\newtheorem*{rem}{Remark}
\newtheorem*{rems}{Remarks}
\numberwithin{equation}{section}

\newcommand{\R}{\mathbb{R}}
\newcommand{\ZZ}{\mathbb{Z}}
\newcommand{\Lor}{\R^{m+2,2}}
\newcommand{\Euc}{\R^{m+2}}
\newcommand{\Q}{\mathcal{Q}}
\newcommand{\Z}{\mathcal{Z}}
\newcommand{\LC}{\mathcal{L}}
\newcommand{\sph}{S^{m+1}}
\newcommand{\mfd}{M^m}
\newcommand{\smfd}{\Sigma^k}
\newcommand{\nml}{\mathcal{N}}
\newcommand{\V}{\mathcal{V}}
\renewcommand{\d}{\mathrm{d}}
\renewcommand{\r}{r} 
\renewcommand{\O}{\mathrm{O}}
\renewcommand{\o}{\mathfrak{o}}
\newcommand{\End}{\mathrm{End}}
\newcommand{\Hom}{\mathrm{Hom}}
\newcommand{\Ah}{\hat{A}}
\newcommand{\fh}{\hat{f}}
\newcommand{\ff}{\mathfrak{f}}
\newcommand{\ffh}{{\hat{\ff}}}  
\newcommand{\sh}{\hat{s}}
\newcommand{\sgh}{\hat{\sigma}}
\newcommand{\hph}{{\hat{\phi}}}
\newcommand{\hn}{\hat{\nu}}
\newcommand{\Vh}{\hat{V}}
\newcommand{\wf}{\Span{\phi,\hph}} 
\newcommand{\set}[1]{\{#1\}}
\newcommand{\Span}[1]{\langle#1\rangle}
\newcommand{\proj}{\mathbb{P}}
\newcommand{\tp}{\Span{t_1}^\perp}
\newcommand{\half}{\tfrac{1}{2}}
\DeclareSymbolFont{script}{U}{eus}{m}{n}
\DeclareSymbolFontAlphabet{\mathscr}{script}
\DeclareMathSymbol{\Wedge}{0}{script}{"5E}
\newcommand{\Rib}{Ribaucour}
\newcommand{\bpt}{Bianchi Permutability Theorem}
\begin{document}
\begin{abstract}
We discuss the Ribaucour transformation of Legendre (contact) maps in its
natural context: Lie sphere geometry.  We give a simple conceptual proof
of Bianchi's original Permutability Theorem and its generalisation by
Dajczer--Tojeiro as well as a higher dimensional version with the
combinatorics of a cube.  We also show how these theorems
descend to the corresponding results for submanifolds in space forms.
\end{abstract}
\maketitle

\section{Introduction}

A persistent and characteristic feature of integrable submanifold
geometries is the existence of transformations of solutions.
Examples include the B\"acklund transformations of surfaces of
constant Gauss curvature and their generalisations
\cite{Bac83,DajToj03,Ten02}; Darboux transformations of isothermic
surfaces \cite{Bia05A,Bur00,Dar99,Sch01}; Eisenhart transformations
of Guichard surfaces \cite[\S 92]{Eis62}; Jonas transformations of
$R$-congruences\footnote{In fact, Jonas transformations are Darboux
transformations of isothermic surfaces in $\R^{2,2}$ \cite{Cal28}.},
to name but a few.

In all these cases, the transformation constructs new submanifolds of
the desired kind from a known one with the help of a solution of an
auxiliary completely integrable first order system of PDE.  Moreover,
some version of the \emph{\bpt} holds:
\begin{quote}
{\itshape Given two transforms of a submanifold there is
a fourth submanifold that is a simultaneous transform of
the first two.}
\end{quote}
This fourth submanifold is often unique and algebraically determined
by the first three.  We say that four submanifolds in this
configuration form a \emph{Bianchi quadrilateral}.

By the 1920's, it was realised that all these transformations had a
common geometric foundation: they were all either \emph{Ribaucour
transformations} or \emph{$W$-transformations} which latter are the
projective-geometric analogue of the former under Lie's line-sphere
correspondence (cf.~\cite{Eis62}).  The \Rib\ transformation was
investigated in classical times for surfaces
\cite{Bia23,Bla29,Eis62}, and for triply and $n$-ply orthogonal
systems \cite{Dem10,Bia23}. Moreover, the integrable systems approach
to orthogonal systems and to discrete orthogonal (circular) nets has
led to renewed interest in the \Rib\ transformation in modern times
\cite{BobHer01,DolSan99,DolManSan98,GanTsa96,KonSch98,Kri97,Zak98}.

It is therefore the purpose of this paper to give a modern treatment
of \Rib\ transformations and their permutability. Let us begin by
describing what they are.  Contemplate a pair of immersions
$\ff,\ffh:\smfd\to N$ of a $k$-manifold into a space-form.  We say
that $\ffh$ is a \emph{\Rib\ transform} of $\ff$ if
\begin{enumerate}
\item[(i)] for each $p\in\smfd$, there is a $k$-sphere $S(p)$ having
first order contact with both $\ff$ and $\ffh$ at $p$.
\item[(ii)] the shape operators of $\ff$ and $\ffh$ commute.
\end{enumerate}
Thus, in classical language, $\ff$ and $\ffh$ parametrise 
two submanifolds of a space form which envelop a congruence of spheres in
such a way that curvature directions at corresponding points
coincide.

We emphasise that these conditions are relatively mild: any
submanifold is enveloped by sphere congruences (in Euclidean space,
the congruence of (one-point compactified) tangent spaces is an
example) and, at least in codimension one, any sphere congruence
generically envelops two (possibly complex) submanifolds.  Moreover,
the condition on curvature directions is not difficult to arrange:
any parallel submanifold to a given submanifold is a \Rib\ transform
(in higher codimension, this example requires that the normal vector
field joining the submanifolds be covariant constant).  Thus the
geometry of \Rib\ transforms is much less rigid than the integrable
specialisations mentioned above.  This is reflected in the following
more general version of the permutability theorem which was proved by
Bianchi \cite[\S 354]{Bia23}:
\begin{BPT}
Given two Ribaucour transforms of a surface there are,
generically, two $1$-parameter families, \emph{Demoulin families}, of
surfaces, one containing the original surface and the other
containing its two Ribaucour transforms, so that any member of one
family is a \Rib\ transformation of any member of the other.

Moreover, corresponding points on any of the surfaces in either
family are concircular.
\end{BPT}
This result has recently been generalised to submanifolds in
space-forms of arbitrary signature with arbitrary dimension and
co-dimension by Dajczer--Tojeiro \cite{DajToj02,DajToj03}.

Moreover, in the context of triply orthogonal systems, a higher
dimensional version of the Bianchi Permutability Theorem has been
obtained by Ganzha--Tsarev \cite{GanTsa96}:
\begin{CT}
Given three initial Ribaucour transforms of a triply orthogonal
system, a generic choice of three simultaneous Ribaucour transforms of
two of them leads to an eighth orthogonal system which is a
simultaneous Ribaucour transform of the latter three thus yielding the
combinatorics of a cube --- a ``Bianchi cube''.
\end{CT}
This higher dimensional version then generalises to any number of
initial Ribaucour transforms and, in this way, gives rise to discrete
orthogonal nets of any dimension by repeatedly applying the
permutability theorem.
In fact, this type of permutability theorem is central to integrable
discretizations of smooth (geometric) integrable systems: it amounts
to the ``consistency condition'' at the heart of the beautiful theory
of Bobenko--Suris
\cite{BobSur02,BobSur05}.
Conversely, this theorem for discrete nets
can be used to prove the permutability theorem for orthogonal systems
of holonomic submanifolds by taking an appropriate limit
\cite{BobMatSur03}.

The aim of this paper is to give a transparent and almost elementary
proof of these permutability theorems for submanifolds in the realm
of Lie sphere geometry.  This is the correct context for a discussion
of the Ribaucour transformation since the main ingredients of the
theory, sphere congruences and the curvature directions of enveloping
submanifolds, are both Lie invariant notions.  Our results descend to
the respective permutability theorems for submanifolds in Riemannian
space-forms (when one imposes the obvious regularity assumptions) and
the Lie invariance of the constructions becomes manifest.

Here is how we will proceed: in Section \ref{sec:1} we swiftly
rehearse the basics of Lie sphere geometry.  Thus a sphere congruence
is viewed as a map into a certain quadric and contact lifts of
submanifolds as Legendre maps into the space of lines in that
quadric.  The enveloping relation is now one of incidence while the
\Rib\ condition amounts to flatness of a certain normal bundle.  In
Section~\ref{sec:perm-theor}, we first prove the Bianchi
Permutability Theorem for Legendre maps which we see amounts to the
assertion that a certain rank four bundle is flat.  Then we prove the
analogue of the Cube Theorem of Ganzha--Tsarev for Legendre maps.  At
this point, our Lie sphere geometric discussion of the Permutability
Theorems is complete: it only remains to show how our results imply
those for submanifolds of space-forms.  To this end, we briefly
discuss \Rib\ transforms of submanifolds of a sphere in
Section~\ref{sec:rib-transf-riem} and then, in
Section~\ref{sec:from-riemannian-lie}, show how to construct Legendre
maps (from the unit normal bundles of our submanifolds) and so find
ourselves in the context of our main results.  Here we make the
effort to work in arbitrary codimension as applications such as those
of \cite{DajToj02} to submanifolds of constant sectional curvature
require this.  Finally, we conclude with a short discussion of an
example in Section~\ref{sec:example}.

\begin{rem}
For clarity of exposition, we have limited ourselves to the case of
definite signature but, modulo the imposition of additional assumptions
of a generic nature, our entire analysis goes through in arbitrary
signature.  In particular, when applied to the Klein quadric (the
projective light cone of $\R^{3,3}$), whose space of lines is the
space of contact elements of $\mathbb{RP}^3$, we recover Bianchi's
Permutability Theorem for focal surfaces of $W$-congruences.
\end{rem}

\section{\Rib\ sphere congruences}
\label{sec:1}

Contemplate $\Lor$: an $(m+4)$-dimensional vector space with metric
$(\,,\,)$ of signature $(m+2,2)$.  Let $\Q$ denote the projective
light-cone of $\Lor$:
\begin{equation*}
\Q=\proj(\LC)=\set{\R v\subset\Lor\colon (v,v)=0, v\neq0}.
\end{equation*}
Thus $\Q$ is a manifold of dimension $m+2$ with a homogeneous action of
$\O(m+2,2)$.

Further, let $\Z$ denote the set of lines in $\Q$ or, equivalently, the
Grassmannian of null $2$-planes in $\Lor$.  Then $\Z$ is an
$\O(m+2,2)$-homogeneous contact manifold of dimension $2m+1$.

The viewpoint of Lie sphere geometry is that $\Q$ parametrises
(non-canonically) the set of oriented hyperspheres in $\sph$
including those of zero radius.  In this picture,
$v,w\in\LC\setminus\set{0}$ are orthogonal if and only if the
corresponding hyperspheres are in oriented contact and so $\Z$
parametrises the space of contact elements on $\sph$ (this is the
origin of the contact structure of $\Z$).

Accordingly, we are lead to the following definitions:
\begin{defn}
A \emph{sphere congruence} is a map $s:\mfd\to\Q$ of an $m$-manifold.
\end{defn}
\begin{defn}
A map $f:\mfd\to\Z$ of an $m$-manifold is a \emph{Legendre map} if,
for all $\sigma_0,\sigma_1\in\Gamma f$, 
\begin{equation*}
(\d \sigma_0,\sigma_1)\equiv 0.
\end{equation*}
\end{defn}
Here, and below, we identify a map $f$ of $\mfd$ into a (subspace of a)
Grassmannian with the corresponding subbundle, also called $f$ of
$\mfd\times\Lor$.  The space of sections $\Gamma f$ therefore
consists of maps $\sigma:\mfd\to\Lor$ with each $\sigma(p)\in f(p)$,
for $p\in\mfd$.

The Legendre condition asserts that $\d f(T\mfd)$ lies in the contact
distribution along $f$ and has the interpretation that $f$ is the
contact lift of a (tube around a) submanifold of $\sph$ (see
Section~\ref{sec:from-riemannian-lie} below).  This motivates our
next definition:
\begin{defn}
A Legendre map $f:\mfd\to\Z$ \emph{envelops} a sphere congruence
$s:\mfd\to\Q$ if $s(p)\subset f(p)$, for all $p\in\mfd$.
\end{defn}

We are interested in the situation where two pointwise distinct
Legendre maps $f,\fh:\mfd\to\Z$ envelop a common sphere congruence
$s:\mfd\to\Q$.  Thus $s=f\cap\fh$. Set
\begin{equation*}
\nml_{f,\fh}=(f+\fh)/s.
\end{equation*}
Thus $\nml_{f,\fh}$ is a rank $2$ subbundle of $s^\perp/s$.  Now
$s^\perp/s$ inherits a metric of signature $(m+1,1)$ from that of
$\Lor$ and this metric restricts to one of signature $(1,1)$ on
$\nml_{f,\fh}$ (otherwise each $f(p)+\fh(p)$ would be a null $3$-plane
in $\Lor$).  In particular, we have a well-defined orthogonal
projection $\pi:s^\perp/s\to\nml_{f,\fh}$.  Observe that, for
$\sigma\in\Gamma s$, the Legendre condition on $f,\fh$ gives
\begin{equation*}
\d\sigma\perp f+\fh.
\end{equation*}
From this we see that, for $\nu\in\Gamma(f+\fh)$,
\begin{equation*}
(\d\nu,\sigma)=-(\nu,\d\sigma)=0
\end{equation*}
so that $\d\nu$ takes values in $s^\perp$ while $\pi(\d\sigma+s)=0$.
We therefore conclude:
\begin{lemma}
There is a metric connection $\nabla^{f,\fh}$ on $\nml_{f,\fh}$ given
by
\begin{equation*}
\nabla^{f,\fh}(\nu+s)=\pi(\d\nu+s).
\end{equation*}
\end{lemma}

\begin{rem}
In case that $s:\mfd\to\Q$ is an immersion, $\nml_{f,\fh}$ can be
identified with the (weightless) normal bundle of $s$ (with respect
to the $\O(m+2,2)$-invariant conformal structure of signature
$(m+1,1)$ on $\Q$) and we have an amusing identification of
enveloping manifolds of $s$ with null normal lines to $s$.  In this
setting, $\nabla^{f,\fh}$ is the (conformally invariant) normal
connection of the weightless normal bundle \cite{BurCal}.
\end{rem}

We are now in a position to make our basic definition:
\begin{defn}
Given Legendre maps $f,\fh:\mfd\to\Z$ enveloping a sphere congruence
$s:\mfd\to\Q$, that is, $s=f\cap\fh$, we say that $s$ is a
\emph{\Rib\ sphere congruence} if $\nabla^{f,\fh}$ is flat.

In this case, $f,\fh$ are said to be \emph{\Rib\ transforms} of each
other and that $(f,\fh)$ are a \emph{\Rib\ pair}.
\end{defn}

\begin{rem}
As we shall see in Section~\ref{sec:from-riemannian-lie}, this
definition generically amounts to the classical notion of two
$k$-dimensional submanifolds of $\sph$ enveloping a congruence of
$k$-spheres in such a way that curvature directions of corresponding
normals coincide, cf. \cite{Bla29,DajToj02,DajToj03}.
\end{rem}

For future reference, we note:
\begin{lemma}
\label{thm:1}
For $\nu\in\Gamma f$, $\nabla^{f,\fh}(\nu+s)=0$ if and only if
$\d\nu\perp\fh$.
\end{lemma}
\begin{proof}
$\nabla^{f,\fh}(\nu+s)=0$ if and only if $\d\nu +s\perp (f+\fh/s)$
or, equivalently, $\d\nu\perp f+\fh$.  But $\d\nu\perp f$ already
since $f$ is Legendre.
\end{proof}

\section{The Permutability Theorem}
\label{sec:perm-theor}

Let $\fh_0,\fh_1$ be \Rib\ transforms of a Legendre map
$f_0:\mfd\to\Z$.  The familiar assertion of Bianchi permutability is
that there should be another Legendre map $f_1$ which is again a
simultaneous \Rib\ transform of the $\fh_i$.
\begin{figure}[h]
\centering
\includegraphics{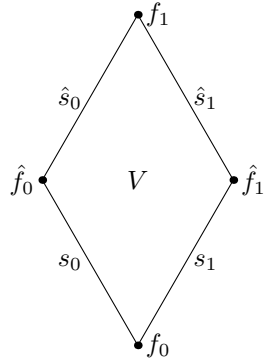}
\caption{A Bianchi quadrilateral}\label{fig:1}
\end{figure}
We say that the $f_i$'s and the $\fh_i$'s form a \textit{Bianchi quadrilateral}
the edges of which represent the enveloped Ribaucour sphere
congruences, see Figure~\ref{fig:1}.
In this section, we characterise the circumstances under which this
assertion holds and show that, in this happy situation, a much
stronger statement is available, cf.~\cite{Bia23}.

For this, we need some brief preliminaries.  So let $(f_0,\fh_0)$,
$(f_0, \fh_1)$ be distinct \Rib\ pairs with sphere congruences $s_i=f_0\cap
\fh_i$.  We impose the mild assumption that $s_0$ and $s_1$ are
pointwise distinct so that $f_0=s_0\oplus s_1$.  It follows that, for
all $p\in\mfd$, $\fh_0\cap\fh_1=\set{0}$ (otherwise
$(\fh_0(p)\cap\fh_1(p))\oplus f_0(p)$ is a null $3$-plane) so that we
may define a rank $4$ subbundle $V$ of $\mfd\times\Lor$ by
\begin{equation*}
V=\fh_0\oplus\fh_1.
\end{equation*}
Since each $\fh_i(p)$ is maximal isotropic in $\Lor$, $V$ inherits a
metric of signature $(2,2)$ from $\Lor$ and thus a metric connection
$\nabla$ by orthoprojection of $\d$.

We may view any section $\sigma_0\in\Gamma s_0$ as representing a
section of $\nml_{f_0,\fh_1}$.  We have:
\begin{lemma}
\label{thm:2}
$\nabla^{f_0,\fh_1}(\sigma_0+s_1)=0$ if and only if $\nabla\sigma_0=0$.
\end{lemma}
\begin{proof}
From Lemma~\ref{thm:1}, we know that
$\nabla^{f_0,\fh_1}(\sigma_0+s_1)=0$ if and only if
$\d\sigma_0\perp\fh_1$.  However, since $\fh_0$ is Legendre, we
already have $\d\sigma_0\perp\fh_0$ whence $\sigma_0+s_1$ is parallel
in $\nml_{f_0,\fh_1}$ if and only if
$\d\sigma_0\perp\fh_0\oplus\fh_1=V$, that is, $\nabla\sigma_0=0$.
\end{proof}

By hypothesis, $\nml_{f_0,\fh_0}$, $\nml_{f_0,\fh_1}$ are both flat
so that we have $\nabla$-parallel sections $\sigma_i\in\Gamma s_i$,
for $i=0,1$.  From this we readily conclude:
\begin{thm}
\label{th:1}
Let $f_1$ be a simultaneous \Rib\ transform of $\fh_0$ and $\fh_1$
which is pointwise distinct from $f_0$.  Then $\nabla$ is flat and
all four subbundles $f_i$, $\fh_i$, $i=0,1$, are $\nabla$-parallel.
\end{thm}
\begin{proof}
Applying the above analysis to $f_1$ in place of $f_0$ provides us
with $\nabla$-parallel sections $\sgh_i\in\Gamma\fh_i$ such that
$f_1=\Span{\sgh_0,\sgh_1}$.  These together with $\sigma_0,\sigma_1$
form a $\nabla$-parallel frame for $V$ and the result follows.
\end{proof}

Thus, if Bianchi permutability holds, $V$ is flat.  Locally, a
converse, indeed, much stronger statement is available: assume that
$V$ is flat and $\mfd$ is simply-connected.  Denote by $\V$ the
vector space of $\nabla$-parallel sections of $V$: evaluation at any
fixed $p\in \mfd$ gives us an isomorphism $\V\cong\R^{2,2}$.  The
projective light cone of $\V$ is a $(1,1)$-quadric and so ruled by
two families of (real) lines (the $\alpha$-lines and $\beta$-lines
of twistor theory).  Each family is parametrised by an $\R P^1$;
lines of the same family do not intersect while each pair of lines
from different families intersects in a unique point (thus our
quadric is an isomorph of $\R P^1\times \R P^1$).

A line in this quadric is the same as a map $f:\mfd\to\Z$ with
$f\subset V$ a $\nabla$-parallel subbundle.
\begin{lemma}
\label{th:2}
Let $f:\mfd\to\Z$ be a map with $f\subset V$,
where $V$ is a flat $(2,2)$-bundle.
Then $f$ is $\nabla$-parallel if and only if $f$ is Legendre.
\end{lemma}
\begin{proof}
Let $\sigma\in\Gamma f$.  Then $\d\sigma\perp f$ if and only if
$\d\sigma$ takes values in $f\oplus V^\perp$ or, equivalently,
$\nabla\sigma$ takes values in $f$. 
\end{proof}

We therefore have:
\begin{thm}
\label{th:3}
If $V$ is a flat $(2,2)$-bundle, then there are two families
$f_\alpha,\fh_\beta:\mfd\to\Z$, $\alpha,\beta\in\R P^1$, of
Legendre maps such that:
\begin{enumerate}
\item[(i)] $f_0,f_1\in\set{f_\alpha}$;
\item[(ii)] $\fh_0,\fh_1\in\set{\fh_\beta}$;
\item[(iii)] $(f_\alpha,\fh_\beta)$ is a \Rib\ pair for all
$\alpha,\beta\in\R P^1$.
\end{enumerate}
\end{thm}
\begin{proof}
Only the third assertion requires any explanation: for this, note that
$\nml_{f_\alpha,\fh_\beta}$ is spanned by sections which are represented
by $\nabla$-parallel sections of $f_\alpha$ and $\fh_\beta$.
By Lemma~\ref{thm:2}, these latter sections are
$\nabla^{f_\alpha,\fh_\beta}$-parallel so that
$\nml_{f_\alpha,\fh_\beta}$ is flat whence
$(f_\alpha,\fh_\beta)$ is a \Rib\ pair.
\end{proof}

We call $\set{f_\alpha}$, $\set{\fh_\beta}$ the \emph{Demoulin
families} of Legendre maps after their discoverer \cite{Dem10}.

It remains to see when this beautiful state of affairs actually occurs:
that is, starting from two Ribaucour transforms $\hat f_i$, $i=0,1$, of a
Legendre map $f_0$, when is $V=\hat f_0\oplus\hat f_1$ flat?
For this, note that $\nabla$ is metric and the sections $\sigma_0,\sigma_1$
are already  $\nabla$-parallel so that there is at most one
direction in $\o(2,2)$ for $R^\nabla$ to take values.  Specifically:

\begin{prop}
\label{th:4}
Let $\sgh_i\in\Gamma \fh_i$ represent spanning sections of
$\fh_i/s_i$, $i=0,1$.  Then $V$ is flat if and only if
$(R^\nabla\sgh_0,\sgh_1)\equiv0$.
\end{prop}

In particular, we can weaken the hypotheses of Theorem~\ref{th:1}:
\begin{prop}
\label{th:8}
$V$ is flat if and only if $f$ admits a Legendre complement
$f_1:\mfd\to\Z$ in $V$: $V=f\oplus f_1$
\end{prop}
\begin{proof}
We can choose sections $\sgh_i$ of $f_1\cap\fh_i$ dual to $\sigma_i$:
$(\sgh_0,\sigma_1)=(\sgh_1,\sigma_0)\equiv1$.  Then
$\d\sgh_0\perp\sigma_0,\sgh_0$ since $\fh_0$ is Legendre, while
\begin{equation*}
(\d\sgh_0,\sigma_1)
   = -(\sgh_0,\d\sigma_1)
   = 0
\end{equation*}
since $\sigma_1$ is $\nabla$-parallel.  Finally $\d\sgh_0\perp\sgh_1$
as $f_1$ is Legendre.  Therefore $\sgh_0$ is $\nabla$-parallel.
Similarly, $\sgh_1$ is parallel so that $V$ is flat since it is
spanned by parallel sections.
\end{proof}

In Section~\ref{sec:from-riemannian-lie}, we shall see that,
generically, the bundle $V$ defined by two Ribaucour transforms of a
Legendre map is automatically flat so that the Permutability Theorem
described in Theorem~\ref{th:3} does indeed hold.  However, we shall
show by an example that $V$ can fail to be flat and so the
Permutability Theorem fails also.

\begin{figure}[ht]
\centering
\includegraphics{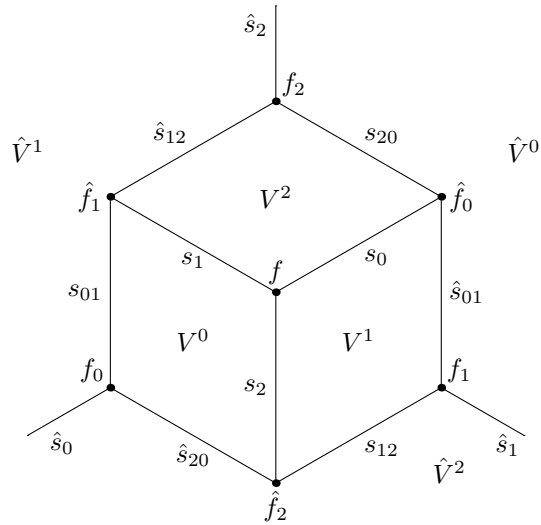}
\caption{A Bianchi cube}
\label{fig:2}
\end{figure}

We now turn to a higher dimensional of the Permutability Theorem
analogous to that obtained in \cite{GanTsa96} for orthogonal systems.
For this, start with a Legendre map $f$ and three Ribaucour
transforms thereof $\fh_0,\fh_1,\fh_2$.  Assume that the Bianchi
Permutability Theorem applies so that we have three more Legendre
maps $f_0,f_1,f_2$ forming three Bianchi quadrilaterals with vertices
$f,\fh_i,f_j,\fh_k$ ($i,j,k$ distinct).  One can now attempt to apply
the theorem again starting with each $\fh_i$ and its Ribaucour
transforms $f_j,f_k$.  The astonishing fact is that there is a single
Legendre map $\fh$ which is a simultaneous \Rib\ transform of all the
$f_i$ so that we obtain a configuration of eight Legendre maps
forming six Bianchi quadrilaterals with the combinatorics of a cube,
a \emph{Bianchi cube}, whose vertices are Legendre maps and whose
edges are the enveloped \Rib\ sphere congruences.  The situation is
illustrated in Figure~\ref{fig:2} where the eighth surface $\fh$ has
been placed at infinity.

This result needs some mild genericity hypotheses.  Here is the
precise statement:
\begin{thm}\label{th:9}
Let $f,\fh_0,\fh_1,\fh_2,f_0,f_1,f_2$ be Legendre maps with each
$f,\fh_i,f_j,\fh_k$ a Bianchi quadrilateral for $i,j,k$ distinct.

Assume that, for $i,j,k$ distinct,
\begin{subequations}
\begin{gather}
\label{eq:7}
\fh_i\not\subset\fh_j\oplus\fh_k\\
\label{eq:8}
f_i\cap f_j=\set{0}\\
\label{eq:9}
f_i\not\subset f_j\oplus f_k.
\end{gather}
\end{subequations}
Then there is a unique Legendre map $\fh$ which is a simultaneous
\Rib\ transform of $f_0,f_1,f_2$.
\end{thm}

\begin{rems}\item{}
\begin{enumerate}
\item[1.] The hypothesis \eqref{eq:7} fails exactly when all the $\fh_i$
lie in a single Demoulin family.  In this degenerate case, any other
$\fh$ in the same family satisfies the conclusion of the theorem
(without the uniqueness assertion!).
\item[2.] The hypothesis \eqref{eq:8} is satisfied for generic choices of
the $f_i$ in their respective Demoulin families.  Indeed, for given
$f_0$, the set of $f_1,f_2$ satisfying \eqref{eq:8} with $f_0$ is
non-empty\footnote{For example, we could take $f_1=f_2=f$} and (Zariski) open.
\item[3.] With a little effort, one can show that hypothesis
\eqref{eq:9} follows from \eqref{eq:7}.
\end{enumerate}
\end{rems}

We begin the proof of Theorem~\ref{th:9} by setting up notation:
given the seven Legendre maps of the statement of the theorem, define
bundles of $(2,2)$-planes by
\begin{align*}
   V^0&:=\fh_1\oplus\fh_2=f\oplus f_0, &
   V^1&:=\fh_2\oplus\fh_0=f\oplus f_1, &
   V^2&:=\fh_0\oplus\fh_1=f\oplus f_2.
\end{align*}
By \eqref{eq:7}, $\fh_i\not\subset V^i$ from which it follows that
\begin{equation*}
V:=\fh_0+\fh_1+\fh_2
\end{equation*}
is a bundle of $(3,2)$-planes\footnote{To see that $V$ has
non-degenerate metric note that any element of $V\cap V^\perp$
together with any $f_i$ would span a null $3$-plane.}.  Further, by
\eqref{eq:8}, we have three more bundles of $(2,2)$-planes:
\begin{align*}
   \Vh^2&:=f_0\oplus f_1,&
   \Vh^0&:=f_1\oplus f_2,&
   \Vh^1&:=f_2\oplus f_0.
\end{align*}
We label the sphere congruences implementing the various Ribaucour
transformations as in Figure~\ref{fig:2}.

Now for the uniqueness assertion: if $\fh$ is a simultaneous \Rib\
transformation of the $f_i$, set $\sh_i=\fh\cap f_i\subset \Vh^i\cap
f_i$.  By \eqref{eq:9}, $f_i\not\subset\Vh^i$ whence $\sh_i=\Vh^i\cap
f_i$ and so is determined by the first seven Legendre maps.  Thus
$\fh$ is so determined also.

To see that a simultaneous \Rib\ transform $\fh$ exists, start by
defining null line-bundles $\sh_i=\Vh^i\cap f_i$.  One readily checks
that, for distinct $i,j,k$, $\sh_i\neq\sh_j$ and
\begin{equation*}
\sh_i\subset \sh_j\oplus\sh_k
\end{equation*}
so that we have a well-defined bundle of $2$-planes
\begin{equation*}
\fh=\sh_i\oplus\sh_j
\end{equation*}
which is null since it contains three distinct null lines (the
$\sh_k$).

In view of Proposition~\ref{th:8}, the only thing left to prove is
that $\fh$ is Legendre.  Our proof of this hinges on the existence of
sections of the various sphere congruences with the property that
their sum around any vertex of the cube vanishes.  \emph{A priori},
this requirement seems overdetermined (even around a single face!)
but we will be able to construct such sections from a consistent
choice of normals to the faces.  This is an entirely algebraic matter
so we consider the situation at a single point.

Thus we contemplate a configuration of null $2$-planes, null lines
and $(2,2)$ planes in $\R^{3,2}$ assigned, respectively, to the
vertices, edges and faces of a combinatorial cube with the property
that each line, resp.\ $(2,2)$-plane, corresponding to an edge,
resp.\ face, is given by the intersection, resp.\ span, of the null
$2$-planes corresponding to incident vertices.  We label the
components of this configuration as in Figure~\ref{fig:2}.

In line with the hypotheses of Theorem~\ref{th:9}, we assume that the
$V^i$ are pairwise distinct, from which it follows that
$V^i\neq\Vh^j$, for all $i,j$, and that the $\Vh^j$ are pairwise
distinct also.

Now let $\nu_i,\hn_j$ be unit normals in $\R^{3,2}$ to $V^i,\Vh^j$.
The situation along an edge of our cube is given by:
\begin{lemma}\label{th:10}
Let $f,\fh_2$ be null $2$-planes and $V^0,V^1$ $(2,2)$-planes in
$\R^{3,2}$ such that
\begin{enumerate}
\item[(i)] $V^0\neq V^1$;
\item[(ii)] $f+\fh_2\subset V^0\cap V^1$;
\item[(iii)] $s_2:=f\cap\fh_2$ is a null line.
\end{enumerate}
Let $\nu_0,\nu_1$ be unit normal vectors to $V^0,V^1$.
Then
\begin{equation*}
\nu_1 = \varepsilon\nu_0 + \sigma,
\end{equation*}
where $\varepsilon=(\nu_0,\nu_1)=\pm1$ and $\sigma_2\in s_2\setminus\{0\}$.
\end{lemma}
\begin{proof}
Write $\nu_1=\varepsilon\nu_0+\sigma_2$ with $\sigma_2\in V^0$ and
$\varepsilon=(\nu_0,\nu_1)$.  Note that $\sigma_2$ is non-zero since
the $\nu_i$ are not collinear.  The $\nu_i$ are orthogonal to
$f,\fh_2$ whence $\sigma_2$ is also.  Thus $\sigma_2$ lies in
$V^0\cap(f+\fh_2)^\perp=s_2$ and we are done since $(\nu_1,\nu_1)=1$
delivers $\varepsilon^2=1$.
\end{proof}

Thus we have a function $\varepsilon$ with values in $\ZZ_2$ defined
on the edges of the cube given by the inner product of adjacent
normals.  Concerning this, we have:
\begin{lemma}\label{th:11}
Let $s_a,s_b,s_c$ be edges meeting at a vertex.  Then
\begin{equation*}
\varepsilon(s_a)\varepsilon(s_b)\varepsilon(s_c)=1.
\end{equation*}
\end{lemma}
\begin{proof}
For definiteness, take the vertex to be $f$.  Apply Lemma~\ref{th:10}
to each edge to get
\begin{align*}
\nu_1&=\varepsilon_2\nu_0+\sigma_2&
\nu_2&=\varepsilon_0\nu_1+\sigma_0&
\nu_0&=\varepsilon_1\nu_2+\sigma_1
\end{align*}
whence
\begin{equation*}
(1-\varepsilon_0\varepsilon_1\varepsilon_2)\nu_1=
\varepsilon_2\varepsilon_1\sigma_0+\varepsilon_2\sigma_1+\sigma_2.
\end{equation*}
The $\sigma_i$ all lie in $f$ so that the right hand side of this is
isotropic and the lemma follows.
\end{proof}

\begin{lemma}
\label{th:12}
There is a choice of normals for which $\varepsilon\equiv1$, that is,
for all $i\neq j$,
\begin{equation*}
(\nu_i,\nu_j)=(\nu_i,\hn_j)=(\hn_i,\hn_j)=1.
\end{equation*}
\end{lemma}
\begin{proof}
Apply Lemma~\ref{th:11} at the vertex $f$: either all $\varepsilon_i$
are $1$ or two of them, $\varepsilon_1,\varepsilon_2$ say, are $-1$
in which case replace $\nu_0$ by $-\nu_0$ to get a choice with all
$\varepsilon_i=1$.

The same argument at $\fh$ gives us $\hn_i$ with all $\varepsilon=1$
around $\fh$.  Moreover, changing the signs of all $\hn_i$ at once,
if necessary, we may assume that all $\varepsilon=1$ around $f_1$.
At each remaining vertex in turn, we have that two of the
$\varepsilon=1$ whence, by Lemma~\ref{th:11}, the third is also.
\end{proof}

\begin{rem}
This argument is cohomological: view $\varepsilon$ as a $1$-cochain on
the octahedron dual to the cube.  Then Lemma~\ref{th:11} says that
$\varepsilon$ is a cocycle while Lemma~\ref{th:12} says that it is a
coboundary (of the $0$-cochain that consistently orients our normals). 
\end{rem}

\begin{prop}\label{th:13}
In the situation of theorem~\ref{th:9}, there are non-zero sections
$\sigma_i,\sgh_j,\sigma_{ij},\sgh_{jk}$ of the participating sphere
congruences the sum of which around any vertex is zero:
\begin{align*}
\sigma_0+\sigma_1+\sigma_2&=0&
\sgh_0+\sgh_1+\sgh_2&=0&
\sigma_{ki}+\sigma_i+\sgh_{ij}&=0&
\sgh_{ki}+\sgh_i+\sigma_{ij}&=0,
\end{align*}
for all $(i,j,k)$ a cyclic permutation of $(1,2,3)$.
\end{prop}
\begin{proof}
With normals chosen as in Lemma~\ref{th:12}, use Lemma~\ref{th:10} to
define sections by:
\begin{align*}
   \nu_0 &= \nu_2 + \sigma_1 = \hat\nu_2 + \hat\sigma_{20}, &
   \hat\nu_0 &= \hat\nu_2 + \hat\sigma_1 = \nu_2 + \sigma_{20}, \\
   \nu_1 &= \nu_0 + \sigma_2 = \hat\nu_0 + \hat\sigma_{01}, &
   \hat\nu_1 &= \hat\nu_0 + \hat\sigma_2 = \nu_0 + \sigma_{01}, \\
   \nu_2 &= \nu_1 + \sigma_0 = \hat\nu_1 + \hat\sigma_{12}, &
   \hat\nu_2 &= \hat\nu_1 + \hat\sigma_0 = \nu_1 + \sigma_{12}.
\end{align*}
These clearly have the desired property.
\end{proof}

With this preparation in hand, we can now complete the proof of
Theorem~\ref{th:9} by showing that $\fh$ is Legendre.  We compute:
\begin{subequations}\label{eq:13}
\begin{align}
\label{eq:10}
(\d\sgh_0,\sgh_1)&=(\d\sigma_{01}+\d\sgh_{20},\sigma_{12}+\sgh_{01})\\
\label{eq:11}
&=(\d\sigma_{01},\sigma_{12})+(\d\sigma_{01},\sgh_{01})+(\d\sgh_{20},\sgh_{01})
\end{align}
\end{subequations}
where we have used Proposition~\ref{th:13} at $f_0,f_1$ for
\eqref{eq:10} and that $f_2$ is Legendre for \eqref{eq:11}.  Now
\begin{equation*}
(\d\sigma_{01},\sigma_{12})=
-(\d\sigma_{01},\sigma_2+\sgh_{20})=-(\d\sigma_{01},\sigma_2)
\end{equation*}
thanks to Proposition~\ref{th:13} at $\fh_2$ and the fact that $f_0$
is Legendre.  Similarly,
\begin{equation*}
(\d\sgh_{20},\sgh_{01})=
-(\d\sigma_2+\d\sigma_{12},\sgh_{01})=
-(\d\sigma_2,\sgh_{01}).
\end{equation*}
On the other hand, since $f_2,f_0$ are Legendre,
\begin{align*}
0=(\d\sgh_{12},\sigma_{20})&=(\d\sigma_{01}+\d\sigma_1,\sigma_0+\sgh_{01})\\
&=(\d\sigma_{01},\sgh_{01})+(\d\sigma_1,\sgh_{01})+(\d\sigma_{01},\sigma_0)
\end{align*}
and substituting all this back into \eqref{eq:13} gives
\begin{align*}
(\d\sgh_0,\sgh_1)&=-(\d\sigma_{01},\sigma_2)-(\d\sigma_{01},\sigma_0)
-(\d\sigma_1,\sgh_{01})-(\d\sigma_2,\sgh_{01})\\
&=(\d\sigma_{01},\sigma_1)+(\d\sigma_0,\sgh_{01})=0
\end{align*}
because $\fh_0,\fh_1$ are Legendre.

\begin{rem}
Viewing Theorem \ref{th:9} as a $3$-dimensional version of the
permutability theorem, it is natural to enquire as to whether higher
dimensional versions of the result are available.  This is indeed the
case: firstly the $4$-dimensional version is equivalent to the usual
Bianchi permutability theorem for quadrilaterals of discrete
Ribaucour transforms of $2$-dimensional discrete orthogonal nets and
this has been proved in the context of conformal geometry by the
second author \cite[\S\S8.5.8 and 8.5.9]{Her03} using a rather
intricate but elementary argument relying solely on Miguel's theorem.
Thereafter, a simple induction argument using the uniqueness
assertion of Theorem~\ref{th:9} establishes the result in any
higher dimension.
\end{rem}

\section{\Rib\ transforms in Riemannian geometry}
\label{sec:rib-transf-riem}

Let us now make contact with the more familiar Riemannian geometry of
the unit sphere\footnote{All our constructions have manifest
conformal invariance so we could have chosen any other
$(m+1)$-dimensional space form as our Riemannian context.} $\sph\subset\Euc$.

An immersion $\ff:\smfd\to\sph$ of a $k$-manifold envelops a
congruence of $k$-spheres if, for each $q\in\smfd$, there is a
$k$-sphere $S(q)\subset\sph$ such that
\begin{equation*}
\ff(q)\in S(q),\qquad \d\ff(T_q\smfd)=T_{\ff(q)}S(q).
\end{equation*}
Thus $S(q)$ is the intersection of $\sph$ with an affine
$(k+1)$-plane $\ff(q)+W_0(q)$ with $\d\ff(T_q\smfd)\subset W_0(q)$.

A second (pointwise distinct) immersion $\ffh:\smfd\to\sph$ envelops
the same sphere congruence exactly when $\ff(q)+W_0(q)=\ffh(q)+W_0(q)$ and
$\d\ffh(T_q\smfd)\subset W_0(q)$.  Otherwise said:
\begin{lemma}
\label{th:5}
$\ff,\ffh:\smfd\to\sph$ envelop a common sphere congruence if and
only if 
\begin{equation}
\label{eq:1}
\Span{\d\ff(T\smfd),\ffh-\ff}=\Span{\d\ffh(T\smfd),\ffh-\ff}.
\end{equation}
\end{lemma}
In this situation, we may therefore define $\r\in\Gamma\End(T\smfd)$
and a $1$-form $\alpha$ by $\d\ffh=\d\ff\circ\r+\alpha(\ffh-\ff)$ or,
equivalently,
\begin{equation}
\label{eq:2}
\d\ffh-\alpha\ffh=\d\ff\circ\r-\alpha\ff.
\end{equation}
Taking the norm-squared of this last and subtracting
$\alpha\otimes\alpha$ yields
\begin{equation}
\label{eq:3}
(\d\ffh,\d\ffh)=(\d\ff\circ\r,\d\ff\circ\r)
\end{equation}
so that $\r$ intertwines the metrics on $\smfd$ induced by $\ff$ and
$\ffh$.  This tensor has a basic role to play in what follows.

Now denote by $W_0$ the rank $(k+1)$ subbundle of $\smfd\times\Euc$
defined by \eqref{eq:1}. Moreover, let $\rho(q)$ be reflection in the
hyperplane orthogonal to $\ffh(q)-\ff(q)$ so that
$\rho:\smfd\to\O(m+2)$.  Then $\rho\ff=\ffh$ and $\rho W_0=W_0$
(since $\ffh-\ff\in\Gamma W_0$).  Now
$\d\ff(T\smfd)=W_0\cap\Span{f}^\perp$ and similarly for $\ffh$ whence
$\rho\d\ff(T\smfd)=\d\ffh(T\smfd)$ and we conclude that $\rho$
provides a metric isomorphism\footnote{With a little more effort, one
can also show that $\rho$ intertwines the normal connections on
$N_\ff$ and $N_\ffh$ also.}  $\rho:N_\ff\to N_\ffh$ between the
normal bundles of $\ff$ and $\ffh$.  This will allow us to compare
curvature directions of $\ff$ and $\ffh$.

For this, let $\xi\in\Gamma N_\ff$ and set
$\lambda=(\xi,\ffh)/\bigl((\ff,\ffh)-1\bigr)$.  Then
$\xi-\lambda\ff\perp \ffh-\ff$ and so is fixed by $\rho$:
\begin{equation}
\label{eq:4}
\xi-\lambda\ff=\rho\xi-\lambda\ffh.
\end{equation}
Together with \eqref{eq:2}, this yields
\begin{equation*}
\bigl(\d(\rho\xi-\lambda\ffh),\d\ffh-\alpha\ffh\bigr)=
\bigl(\d(\xi-\lambda\ff),\d\ff\circ\r-\alpha\ff\bigr)
\end{equation*}
to both sides of which we add $\d\lambda\otimes\alpha$ to get
\begin{equation*}
\bigl(\d(\rho\xi-\lambda\ffh),\d\ffh\bigr)=
\bigl(\d(\xi-\lambda\ff),\d\ff\circ\r\bigr).
\end{equation*}
We express this last in terms of the shape operators $A$, $\Ah$ of
$\ff$ and $\ffh$:
\begin{equation*}
\bigl(\d\ffh(\Ah^{\rho\xi}+\lambda),\d\ffh\bigr)=
\bigl(\d\ff(A^\xi+\lambda),\d\ff\circ\r\bigr)
\end{equation*}
and use \eqref{eq:3} to conclude that $\r$ also intertwines shape
operators:
\begin{equation}
\label{eq:5}
\r\circ(\Ah^{\rho\xi}+\lambda)=A^{\xi}+\lambda.
\end{equation}
Now $\Ah^{\rho\xi}$ is symmetric with respect to the metric induced
by $\ffh$ so that (computing transposes with respect to the metric
induced by $\ff$)
\begin{equation*}
\r^T\r\Ah^{\rho\xi}=(\Ah^{\rho\xi})^T\r^T\r,
\end{equation*}
which, together with \eqref{eq:5}, yields
\begin{equation}
\label{eq:6}
\r^T(A^\xi+\lambda)=(A^\xi+\lambda)\r.
\end{equation}
With this in hand, we can explain the significance of $\r$ for us:
\begin{prop}
\label{th:6}
If $\r$ is symmetric with respect to the metric induced by $\ff$ then
corresponding shape operators of $\ff$ and $\ffh$ commute:
$[A^\xi,\Ah^{\rho\xi}]=0$, for all $\xi\in N_{\ff}$.

Conversely, if, for some $\xi\in N_{\ff}$, $[A^\xi,\Ah^{\rho\xi}]=0$
and, additionally, $A^\xi+\lambda$ is invertible\footnote{When
$\smfd$ is a hypersurface, this is precisely the condition that our
sphere congruence contains no curvature spheres.}, then $\r$ is symmetric.
\end{prop}
\begin{proof}
From \eqref{eq:5} we have
\begin{equation*}
0=[A^\xi+\lambda,\r](\Ah^{\rho\xi}+\lambda)+\r[A^\xi,\Ah^{\rho\xi}].
\end{equation*}
Now, if $\r$ is symmetric, \eqref{eq:6} gives $[A^\xi+\lambda,\r]=0$
and the conclusion follows from the invertibility of $\r$.

For the converse, given $\xi\in N_{\ff}$ with
$[A^\xi,\Ah^{\rho\xi}]=0$ and $A^\xi+\lambda$, equivalently
$A^\xi+\lambda$, invertible, we first deduce that
$[A^\xi+\lambda,\r]=0$ and then, from \eqref{eq:6}, that $\r$ is symmetric.
\end{proof}

The notion of \Rib\ transform currently available in the literature
\cite{Bla29,DajToj02,DajToj03,Eis62,Ten02} involves a pair of $k$-dimensional
submanifolds enveloping a congruence of $k$-spheres so that curvature
directions of corresponding normals coincide (that is, corresponding
shape operators commute).  In view of Proposition~\ref{th:6}, we
propose the following
\begin{defn}
Immersions $\ff,\ffh:\smfd\to\sph$ enveloping a congruence of
$k$-spheres are a \emph{\Rib\ pair} if $\r$ is symmetric.
\end{defn}

\begin{rem}
Given such a \Rib\ pair, every $A^\xi$ commutes with $\r$.  Thus, in
the generic case where $\r$ has distinct eigenvalues, all shape
operators of $\ff$ must commute with each other so that the normal
bundle of $\ff$ is flat (whence $\ffh$ has flat normal bundle also).
\end{rem}

\section{From Riemannian to Lie sphere geometry }
\label{sec:from-riemannian-lie}

Given $\ff,\ffh:\smfd\to\sph$ enveloping a congruence of $k$-spheres,
we are going to construct Legendre maps $f,\fh:\mfd\to\Z$ enveloping
$s:\mfd\to\Q$ where $\mfd$ is the unit normal bundle of $\ff$.

We shall show:
\begin{enumerate}
\item[(i)] $\nml_{f,\fh}$ is flat if and only if $\r$ is symmetric so that
our two notions of \Rib\ pair correspond.
\item[(ii)] Given \Rib\ transforms $\ffh_0,\ffh_1$ with corresponding
$\r_0,\r_1$, the bundle $V=\fh_0\oplus \fh_1$ of
Section~\ref{sec:perm-theor} is flat if and only if $[\r_0,\r_1]=0$.
\item[(iii)] In this latter situation, all the Legendre maps $f_s,\fh_t$
participating in the Permutability Theorem also arise from maps
$\ff_s,\ffh_t:\smfd\to\sph$.
\end{enumerate}

For all this, fix orthogonal unit time-like vectors $t_0,t_1\in\Lor$
and set $\Euc=\Span{t_0,t_1}^\perp$.  Thus
\begin{equation*}
\Lor=\Euc\oplus\Span{t_0}\oplus\Span{t_1}
\end{equation*}
is an orthogonal decomposition.  The quadric $\Q$ splits as a
disjoint union $\Q=\Q_0\cup\Q_+$ where
\begin{equation*}
\Q_0=\set{\Span{v}\in\Q:v\perp t_1}
\end{equation*}
is the space of \emph{point spheres}.
Note that $x\mapsto\Span{x+t_0}:\sph\to\Q_0$ is a diffeomorphism.  We
therefore define $\phi,\hph:\smfd\to\LC\subset\Lor$ by
\begin{equation*}
\phi=\ff+t_0,\qquad\hph=\ffh+t_0.
\end{equation*}

Now let $\mfd$ be the unit normal bundle of $\ff$ with bundle
projection $\pi:\mfd\to\smfd$.  For $\xi\in\mfd$ with $\pi(\xi)=q$,
define $\lambda(\xi)$ by
\begin{equation*}
\lambda(\xi)=(\xi,\hph(q))/(\phi(q),\hph(q))=
(\xi,\fh(q))/\bigl((\ff(q),\ffh(q))-1\bigr).
\end{equation*}
Observe that $\xi-\lambda(\xi)\phi(q)\perp\phi(q),\hph(q)$ while
\eqref{eq:4} gives:
\begin{equation*}
\xi-\lambda(\xi)\phi(q)=\rho\xi-\lambda(\xi)\hph(q).
\end{equation*}
Here the left hand side is clearly orthogonal to
$\d\phi(T_q\smfd)=\d\ff(T_q\smfd)$ while the right hand side is
orthogonal to $\d\hph(T_q\smfd)$.  Thus, defining
$\sigma:\mfd\to\LC\subset\Lor$ by
\begin{equation*}
\sigma(\xi)=\xi-\lambda(\xi)\phi(q)+t_1=\rho\xi-\lambda(\xi)\hph(q)+t_1,
\end{equation*}
we readily conclude:
\begin{align}
(\phi\circ\pi,\sigma)&=(\hph\circ\pi,\sigma)=0\\
(d(\phi\circ\pi),\sigma)
&=(d(\hph\circ\pi),\sigma)=0
\end{align}
We therefore have Legendre maps $f,\fh:\mfd\to\Z$ with
$f\cap\fh=s:\mfd\to\Q$ as follows:
\begin{equation*}
f=\Span{\phi\circ\pi,\sigma},\qquad
\fh=\Span{\hph\circ\pi,\sigma},\qquad
s=\Span{\sigma}.
\end{equation*}

\begin{rem}
Here is the geometry of the situation: a unit normal $\xi$ to $\ff$
at $q$ defines a contact element $\Span{\xi}^\perp\subset
T_{\ff(q)}\sph$ containing $\d\ff(T_q\smfd)$: this is $f(\xi)$.
Among the hyperspheres sharing this contact element is exactly one
which is also tangent to $\ffh$ at $q$: this is $s(\xi)$.
\end{rem}

We have (at last!) found ourselves in the setting of
Section~\ref{sec:1} and so can investigate $\nml_{f,\fh}$.  For this,
first contemplate the bundle $\wf\to\smfd$ with metric and connection
$\nabla^{\phi,\hph}$ inherited from $\smfd\times\Lor$.  Using
$(\sigma,t_1)\equiv -1$ and the fact that
$\phi\circ\pi+s$, $\hph\circ\pi+s$ are spanning sections of
$\nml_{f,\fh}$, it is not difficult to show that:
\begin{lemma}
There is a metric, connection-preserving isomorphism between
$\nml_{f,\fh}$ and the pull-back $\pi^{-1}\wf$ given by
\begin{equation*}
\tau+s\mapsto \tau+(\tau,t_1)\sigma.
\end{equation*}
\end{lemma}
In particular, since $\pi$ is a submersion, we conclude:
\begin{prop}
$\nml_{f,\fh}$ is flat if and only if $\wf$ is flat\footnote{This
last is the definition of \Rib\ transform adopted by
Burstall--Calderbank \cite{BurCal} in their conformally invariant
treatment of this topic.}
\end{prop}
The latter condition is easy to characterise: define the second
fundamental form $\beta\in\Omega^1_{\smfd}\otimes\Hom(\wf,\wf^\perp)$ of
$\wf$ by
\begin{equation*}
\d\psi=\nabla^{\phi,\hph}\psi+\beta\psi,
\end{equation*}
for $\psi\in\Gamma\wf$, and deduce the following Gauss equation from
the flatness of $\d$:
\begin{equation*}
R^{\nabla^{\phi,\hph}}=\beta^T\wedge\beta.
\end{equation*}
Thus flatness of $\wf$ amounts to the vanishing of
$(\beta\hph\wedge\beta\phi)$.  However, \eqref{eq:2} gives
\begin{equation*}
\d\hph-\alpha\hph=\d\phi\circ\r-\alpha\phi
\end{equation*}
with both sides palpably orthogonal to $\wf$ whence
\begin{equation*}
\beta\hph=\d\hph-\alpha\hph,\qquad\beta\phi=\d\phi-(\alpha\circ\r^{-1})\phi
\end{equation*}
and, in particular, $\beta\hph=(\beta\circ\r)\phi$. Thus
\begin{equation*}
(\beta\hph\wedge\beta\phi)=((\beta\circ\r)\phi\wedge\beta\phi)=
(\d\ff\circ\r\wedge\d\ff)
\end{equation*}
and we conclude:
\begin{thm}
$\nml_{f,\fh}$ is flat if and only if $\r$ is symmetric.
\end{thm}
That is $(f,\fh)$ is a \Rib\ pair of Legendre maps if and only if
$(\ff,\ffh)$ are a \Rib\ pair of submanifolds.

Suppose now that we are in the situation of the permutability
theorem: thus we are given two \Rib\ transforms
$\ffh_0,\ffh_1:\smfd\to\sph$ of $\ff_0:\smfd\to\sph$.  We therefore have
$\r_0,\r_1\in\Gamma\End(T\smfd)$ symmetric with respect to the metric
induced by $\ff_0$.

Assume, once and for all, that these three maps are pairwise
pointwise distinct.  This ensures that $\fh_0\cap\fh_1=\set{0}$ so
that we can define $V=\fh_0\oplus\fh_1$.
\begin{thm}
\label{th:7}
$V$ is flat if and only if $[\r_0,\r_1]=0$.
\end{thm}
\begin{proof}
With $s_i=\fh_i\cap f_0$ and $\hph_i=\ffh_i+t_0$, $\hph_i\circ\pi$
represent non-zero sections of $\fh_i/s_i$ so that, by
Lemma~\ref{th:4}, $V$ is flat if and only if
$\bigl(R^\nabla(\hph_0\circ\pi),\hph_1\circ\pi\bigr)$ vanishes.  Let
$\beta_V\in\Omega_{\mfd}\otimes\Hom(V,V^\perp)$ be the second
fundamental form of $V$ so that
\begin{equation*}
\d\psi=\nabla\psi+\beta_V\psi,
\end{equation*}
for $\psi\in\Gamma V$.  As before, a Gauss equation gives
$R^\nabla=\beta_V^T\wedge\beta_V$ so that flatness of $V$ amounts to
the vanishing of
$\bigl(\beta_V(\hph_0\circ\pi)\wedge\beta_V(\hph_1\circ\pi)\bigr)$.

On the other hand, we also have
\begin{equation*}
\d\hph_i=\nabla^{\phi_0,\hph_i}\hph_i+\beta^i\hph_i
\end{equation*}
with $\beta^i\in\Omega^1_{\smfd}\otimes\Hom(\Span{\phi_0,\hph_i},
\Span{\phi_0,\hph_i}^\perp)$.  Now $\pi^{-1}\Span{\phi_0,\hph_i}\subset
V$ so that $(\pi^*\beta^i-\beta_V)\hph_i\circ\pi$ takes values in
$V$.  Moreover, each $\fh_i$ is Legendre so that
$\d(\hph_i\circ\pi)\perp s_i$ as is $\nabla^{\phi_0,\hph_i}\hph_i$.
We conclude that
$(\pi^*\beta^i-\beta_V)\hph_i\circ\pi\perp\phi_0,\hph_i,s_i$ and so
takes values in $V\cap\Span{\phi_0,\hph_i,s_i}^\perp=s_i$.  Since
$V^\perp, s_0, s_1$ are mutually orthogonal, this gives
\begin{align*}
\bigl(\beta_V(\hph_0\circ\pi)\wedge\beta_V(\hph_1\circ\pi)\bigr)&=
\pi^*\bigl(\beta^0\hph_0\wedge\beta^1\hph_1\bigr)\\
&=\pi^*\bigl((\beta^0\circ\r_0)\phi_0\wedge(\beta^1\circ\r_1)\phi_0\bigr)\\
&=\pi^*\bigl(\d\ff_0\circ\r_0\wedge\d\ff_0\circ\r_1\bigr)
\end{align*}
whence the conclusion.
\end{proof}

Let us suppose then that $[\r_0,\r_1]$ vanishes so that, by Theorems
\ref{th:3} and \ref{th:7}, the Permutability Theorem holds.  It
remains to show that all the (locally defined) Legendre maps
$f_\alpha,\fh_\beta:\mfd\to\Z$ of Theorem~\ref{th:3} arise from maps
$\ff_\alpha,\ffh_\beta:\smfd\to\sph$.  That is, each \emph{point
map} $f_\alpha\cap\Span{t_1}^\perp,\fh_\beta\cap\Span{t_1}^\perp$ is
constant on the fibres of $\pi$.  One approach to this, valid on
simply-connected subsets of $\mfd$ which have connected intersections
with the fibres of $\pi$, is to compute derivatives of these point
maps along said fibres.  However, we employ an alternative, slightly
indirect, argument which provides us with rather more information: we
will show that all point maps in a Demoulin family arise as
pull-backs of sections of a certain bundle of $(2,1)$-planes
in $\tp$ which are parallel with respect to a certain metric
connection.   The bundle is the same for each Demoulin family but the
connections are different.

For all this, we begin by recalling the well-known fact that a
$2$-plane $\Span{\sigma,\tau}\subset V$ is null if and only if the
$2$-vector $\sigma\wedge\tau\in\Wedge^2 V$ is self-dual or
anti-self-dual.  Choose the orientation on $V$ for which $\Wedge^2
f_0$ is self-dual and let $*:\Wedge^2 V\to \Wedge^2 V$ be the
corresponding Hodge star operator with orthogonal eigenspace
decomposition
\begin{equation*}
\Wedge^2 V=\Wedge^2_+V\oplus\Wedge^2_- V.
\end{equation*}
The $(2,2)$-metric on $V$ induces a $(2,4)$-metric from which both
$\Wedge^2_\pm V$ inherit a $(1,2)$-metric.  The flat metric connection
$\nabla$ on $V$ induces flat metric connections $\nabla^\pm$ on
$\Wedge^2_\pm V$ and the isotropic parallel line subundles with
respect to $\nabla^+$, respectively $\nabla^-$, are the $\Wedge^2
f_\alpha$, respectively $\Wedge^2 \fh_\beta$.

Now contemplate the bundle $U=\Span{\phi_0,\hph_0,\hph_1}$: a bundle
of $(2,1)$-planes in $\tp$ over $\smfd$.  Note that
$\pi^{-1}U=V\cap\tp$.  Here is the geometry of $U$: under
the diffeomorphism $x\mapsto\Span{x+t_0}$ of $\sph$ with $\Q_0$,
$U(p)\cap\Q_0$ represents the intersection of $\sph$ with the affine
$2$-plane containing the points $\ff_0(p),\ffh_0(p),\ffh_1(p)$, that
is, the circle containing these points.

Choose $t\in\Gamma V$ with $(t,t)=-1$, $t\perp \pi^{-1}U$ so that
$V=\pi^{-1}U\oplus\Span{t}$.  For $\eta\in\Wedge^2 V$, we note that
the interior product $\iota_t\eta\perp t$ whence
$\iota_t\eta\in\pi^{-1}U$.  We therefore define $T_{\pm}:\Wedge^2_\pm
V\to \pi^{-1}U$ by
\begin{equation*}
T_{\pm}\eta=\sqrt{2}\iota_t\eta.
\end{equation*}
Since $*$ is an involutive isometry that permutes $\Wedge^2
\pi^{-1}U$ with $\pi^{-1}U\wedge\Span{t}$, we have:
\begin{lemma}
\label{th:14}
$T_{\pm}$ is an anti-isometric\footnote{Thus
$(T_{\pm}\eta_1,T_{\pm}\eta_2)=-(\eta_1,\eta_2)$.} isomorphism with
inverse
\begin{equation*}
T_{\pm}^{-1}\phi=\bigl(\phi\wedge t\pm*(\phi\wedge t)\bigr)/\sqrt{2}.
\end{equation*}
Moreover,
\begin{align*}
T_+\Wedge^2 f_\alpha&=f_\alpha\cap\tp & 
T_-\Wedge^2\fh_\beta&=\fh_\beta\cap\tp.
\end{align*}
\end{lemma}

We use $T_\pm$ to induce flat metric connections, also called
$\nabla^\pm$, on $\pi^{-1}U$.  By Lemma~\ref{th:14}, the parallel
isotropic line subbundles of $\pi^{-1}U$ with respect to $\nabla^+$,
respectively $\nabla^-$,are the point maps $f_\alpha\cap\tp$,
respectively $\fh_\beta\cap\tp$.  This means we will be done as soon as
we prove
\begin{prop}
\label{th:15}
There are flat metric connections $D^\pm$ on $U$ such that
\begin{equation*}
\nabla^\pm=\pi^{-1}D^\pm.
\end{equation*}
\end{prop}

For this, we need the following lemma which is surely well-known
(and, in any case, a straight-forward exercise to prove):
\begin{lemma}
\label{th:16}
Let $\pi:\mfd\to\smfd$ be a bundle with connected fibres, $U\to\smfd$
a vector bundle and $\nabla$ a connection on $\pi^{-1}U$.  Then
$\nabla=\pi^{-1}D$ for some connection $D$ on $U$ if and only if, for all
$X\in\ker\d\pi$ and $\phi\in\Gamma U$,
\begin{enumerate}
\item[(i)] $\nabla_X(\phi\circ\pi)=0$;
\item[(ii)] $\iota_X R^\nabla=0$.
\end{enumerate}
\end{lemma}

With this in hand, we compute: for $\psi\in\Gamma\pi^{-1}U$,
\begin{equation}
\label{eq:12}
\begin{split}
\nabla^\pm\psi&=\iota_t\nabla\bigl(\psi\wedge t\pm *(\psi\wedge
t)\bigr)\\
&=\iota_t\bigl(\nabla\psi\wedge t+\psi\wedge\nabla t\pm
*(\nabla\psi\wedge t+\psi\wedge\nabla t)\bigr)\\
&=\pi^{-1}\nabla^U\psi\pm \iota_t*(\psi\wedge\nabla t)
\end{split}
\end{equation}
where $\nabla^U$ is the connection on $U$ induced by $\d$ (whence
$\pi^{-1}\nabla^U\psi$ is the $U$-component of $\nabla\psi$).  Now
let $X\in\ker\d\pi$ and $\phi\in\Gamma U$.  We know that
$\nabla_X(\phi\circ\pi)$ is the $V$-component of $\d_X(\phi\circ\pi)$
and so vanishes.  It follows that $\nabla_X$ preserves $\pi^{-1}U$
whence $\nabla_X t=0$.  Since both $\nabla^\pm$ are flat, Lemma
\ref{th:16} assures us that Proposition~\ref{th:15} holds.

We have therefore arrived at the following situation: on a simply
connected open subset $\Omega\subset\smfd$, we have isotropic line
subbundles $\Span{\phi_\alpha}$, parallel for $D^+$, and
$\Span{\hph_\beta}$, parallel for $D^-$, so that
\begin{align*}
f_\alpha\cap\tp&=\Span{\phi_\alpha\circ\pi}&
\fh_\beta\cap\tp&=\Span{\hph_\beta\circ\pi}.
\end{align*}
Define $\ff_\alpha,\ffh_\beta:\Omega\to\sph$ by
\begin{align*}
\Span{\ff_\alpha+t_0}&=\Span{\phi_\alpha}&
\Span{\ffh_\beta+t_0}&=\Span{\hph_\beta}
\end{align*}
and finally deduce the Bianchi Permutability Theorem for Ribaucour
transforms of maps $\smfd\to\sph$:
\begin{thm}
\label{th:17}
Let $\ff_0:\smfd\to\sph$ be an immersion of a simply connected
manifold with $\ffh_0,\ffh_1$
pointwise distinct Ribaucour transforms satisfying
$[\r_0,\r_1]=0$.  Then, for $\alpha,\beta\in\R\proj^1$, there are
maps $\ff_\alpha,\ffh_\beta:\smfd\to\sph$ such that:
\begin{enumerate}
\item[(i)] $\ff_0\in\set{\ff_\alpha}$; $\ffh_0,\ffh_1\in\set{\ffh_\beta}$;
\item[(ii)] each $\ff_\alpha$ is a Ribaucour transform of each
$\ffh_\beta$;
\item[(iii)] for each $p\in\smfd$, the points $\ff_\alpha(p),\ffh_\beta(p)$,
$\alpha,\beta\in\R\proj^1$, are concircular;
\item[(iv)] any four maps in one Demoulin family, either $\set{\ff_\alpha}$
or $\set{\ffh_\beta}$, have constant cross-ratio.
\end{enumerate}
\end{thm}
\begin{proof}
Only the last point requires further elaboration: any $\ff_\alpha+t_0$
spans a $D^+$-parallel line bundle which admits a parallel section.
The cross-ratio of four such maps can be computed in terms of the
inner products between these parallel sections \cite[\S6.5.4]{Her03}
and so is constant since $D^+$ is metric.
\end{proof}

\begin{rems}\item{}
\begin{enumerate}
\item[1.] Assertion (iii), that corresponding points of the maps of a
Bianchi quadrilateral lie on circles, provides a more direct approach
to Theorem~\ref{th:9} in the context of M\"obius geometry: one can
construct the eighth map in a Bianchi cube via Miguel's theorem.
This gives a link between the Bianchi Permutability Theorem and the
theory of ``discrete orthogonal nets'', or ``circular nets'', see for
example \cite{BobHer01}.
\item[2.] That corresponding points of members of a single Demoulin
family are concircular is due to Demoulin \cite{Dem10} while Bianchi showed
that the circles for the two families coincide \cite[\S 354]{Bia23}.  
Assertion (iv) on cross-ratios is also due to Demoulin \cite{Dem10A}.
\end{enumerate}
\end{rems}

\begin{rem}
Observe that \eqref{eq:12} tells us that $D^\pm$ are of the form
\begin{equation*}
D^\pm=\nabla^U\pm B
\end{equation*}
for some $B\in\Omega^1\otimes\o(U)$.  Moreover, by construction,
$\Span{\phi_0}$ is $D^+$-parallel while $\Span{\hph_0},\Span{\hph_1}$
are $D^-$-parallel.  These properties fix $B$ (and so $D^\pm$)
uniquely: the difference of two such $B$ would be $\o(U)$-valued
while preserving the decomposition
$\Span{\hph_0}\oplus\Span{\phi_0}\oplus\Span{\hph_1}$ and so must
vanish.  This suggests an alternative approach to the Permutability
Theorem entirely in the context of conformal geometry
\end{rem}

\section{Example}\label{sec:example}

We conclude by presenting a very simple configuration of two Ribaucour
transforms $\ffh_0$ and $\ffh_1$ of a surface $\ff_0$ in $S^3$, where the
\bpt\ can fail.  For this, the curvature directions of $\ff_0$ should
be ambiguous while those of $\ffh_0$ and $\ffh_1$ should be
well-defined and different.   Thus $\ff_0$ should parametrise (part
of) a $2$-sphere $s$ and we will take the $\ffh_i$ to be Dupin
cyclides as these are the simplest surfaces in Lie sphere geometry.
Recall that all Dupin cyclides are equivalent in Lie sphere geometry:
they are congruent to a circle (or, equivalently, a torus of
revolution) \cite{Pin81}.

Thus, let $\Q$ be the projective light-cone of $\R^{4,2}$ and fix a
unit time-like $t_1\in\R^{4,2}$ to get a space
\begin{equation*}
{\Q}_0={\Q}\cap t_1^\perp
\end{equation*}
of point spheres.  We can then write $s=\Span{e+t_1}$ with $e$ a
space-like unit vector in $\Span{t_1}^\perp$.

\subsection{The first Ribaucour transform}
\label{sec:first-riba-transf}

We fix two points on $s$:
\begin{equation*}
   p_0,p_\infty\perp t_1,s,
   \quad
   |p_0|^2=|p_\infty|^2=0,
   \quad
   (p_0,p_\infty)=-1,
\end{equation*}
and choose an orthonormal basis $(e_1,e_2)$ for
$\Span{t_1,e,p_0,p_\infty}^\perp$.
Geometrically, $t_1+e_1$ and $t_1+e_2$ define two $2$-spheres that contain
the points $\Span{p_0}$ and $\Span{p_\infty}$.  The circle in which
these spheres intersect is a (degenerate)
Dupin cyclide which we take as our first Ribaucour transform. 
Thus we define a Legendre map
\begin{equation*}
\fh_0=\Span{\hat\kappa_{01},\hat\kappa_{02}}
\end{equation*}
where
\begin{align*}
(u,v)\mapsto
   \hat\kappa_{01}(v)&:=t_1+\cos v\,e_1+\sin v\,e_2\\
(u,v)\mapsto
   \hat\kappa_{02}(u)&:=p_0+u\,e+\frac{u^2}{2}\,p_\infty\perp t_1.
\end{align*}
A section of the corresponding point map is given by
$\hat\phi_0:=\hat\kappa_{02}$. 

We parametrise (the contact lift of) $s$ by  $f_0:=\Span{s,\phi_0}$,
where (a section of) the point map is given by
\begin{equation*}
\phi_0(u,v)
   := p_0+u\,(\cos v\,e_1+\sin v\,e_2)+\frac{u^2}{2}\,p_\infty
   \perp t_1.
\end{equation*}
Then the section
\begin{equation*}
 \sigma_0
   = (\hat\kappa_{01},s)\hat\kappa_{02}
   - (\hat\kappa_{02},s)\hat\kappa_{01}
   \in\Gamma(\fh_0)
\end{equation*}
defines a sphere congruence that is enveloped by both $\fh_0$
and $f_0$ since $\sigma_0(u,v)=u\,s+\phi_0(u,v)$.

Now define
\begin{equation*}
\r_0:=\frac{\partial}{\partial u}\,\d u
   \quad{\rm and}\quad
   \alpha:=\frac{1}{u}\,\d u
\end{equation*}
and note that
\begin{equation*}
\d\hat\phi_0-\alpha\,\hat\phi_0
   = \frac{1}{u}(-p_0+\frac{u^2}{2}p_\infty)\,\d u
   = \d\phi_0\circ\r_0-\alpha\,\phi_0.
\end{equation*}
We see that $\r_0$ is symmetric with respect to the metric $ \d
u^2+u^2\d v^2 $ induced by $\phi_0$ and so conclude that $\fh_0$ is a
Ribaucour transform\footnote{Alternatively, $e+t_1$ represents a
parallel section of $\nml_{f_0\fh_0}$.}  of $f_0$.

\subsection{The second Ribaucour transform}
\label{sec:second-riba-transf}

Now fix two possibly different points on $s$:
$$
   \tilde p_0,\tilde p_\infty\perp t_1,s,
   \quad
   |\tilde p_0|^2=|\tilde p_\infty|^2=0,
   \quad
   (\tilde p_0,\tilde p_\infty)=-1,
$$
and choose an orthonormal basis $(\tilde e_1,\tilde e_2)$ for
$\Span{ t_1,e,\tilde p_0,\tilde p_\infty}^\perp$. Then $s+\tilde p_0$
and $s+\tilde p_\infty$ define two $2$-spheres which touch $s$ at the
points $\Span{\tilde p_0}$ and $\Span{\tilde p_\infty}$
respectively.  These spheres intersect in a circle which we take as
our second Ribaucour transform.

Thus we define a Legendre map by 
\begin{equation*}
\fh_1=\Span{\hat\kappa_{11},\hat\kappa_{12}}
\end{equation*}
where
\begin{align*}
(u,v)\mapsto
   \hat\kappa_{11}(v)&:= 
(e+\tilde p_0+\tilde p_\infty)+\cos v\,\tilde e_1+\sin v\,\tilde
   e_2\perp t_1\\
   (u,v)\mapsto
   \hat\kappa_{12}(u)&:=(1-u+\frac{u^2}{2})\,t_1
    +(1-u)(e+\tilde p_0)
    +(-u+\frac{u^2}{2})(e+\tilde p_\infty).
\end{align*}
Then $\hph_1:=\hat\kappa_{11}$ is a section of the corresponding
point map.

Now we parametrise $s$ by $\tilde f_0:=\Span{ s,\tilde\phi_0}$ with
point map
\begin{equation*}
\tilde\phi_0(u,v)
   := \tilde p_0+u\,(\cos v\,\tilde e_1+\sin v\,\tilde
   e_2)+\frac{u^2}{2}\,\tilde p_\infty
   \perp t_1.
\end{equation*}
Then the section
\begin{equation*}
   \sigma_1
   = (\hat\kappa_{11},s)\hat\kappa_{12}
   - (\hat\kappa_{12},s)\hat\kappa_{11}
   \in\Gamma(\fh_1)
\end{equation*}
defines a sphere congruence that is enveloped by $\fh_1$ and $\tilde
f_0$ since $\sigma_1(u,v)=(1-u+\frac{u^2}{2})\,s+\tilde\phi_0(u,v)$.

Now define
\begin{equation*}
   \r_1:=\frac{1}{u}\frac{\partial}{\partial v}\,\d v
   \quad{\rm and}\quad
   \alpha:=0
\end{equation*}
and note that
\begin{equation*}
   \d\hat\phi_1-\alpha\,\hat\phi_1
   = (-\sin v\,\tilde e_1+\cos v\,\tilde e_2)\,\d v
   = \d\tilde\phi_0\circ\r_1-\alpha\,\tilde\phi_0.
\end{equation*}
Again, $\r_1$ is symmetric with respect to the metric
$
\d u^2+u^2\d v^2
$
induced by $\tilde\phi_0$ showing that $\fh_0$ is a Ribaucour transform
of $\tilde f_0$.
\begin{figure}[ht]
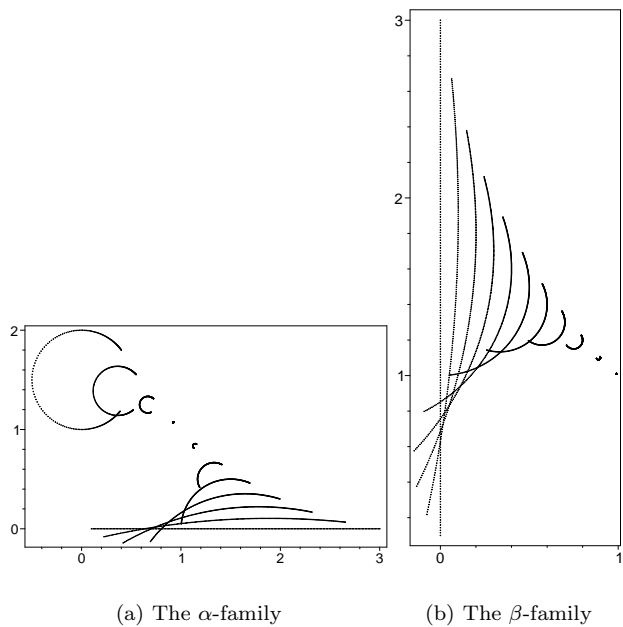

\centering
\subfigure[The $\alpha$-family]{\includegraphics[width=5cm]{graphic-3.eps}}
\subfigure[The $\beta$-family]{\includegraphics[width=3cm]{graphic-4.eps}}
\caption{Demoulin families of Dupin cyclides}\label{fig:3}
\end{figure}

\subsection{The Permutability Theorem}
\label{sec:perm-theor-1}
Now let us see when the \bpt\ holds.  For this we should choose a
common parametrisation of all participating surfaces but we can avoid
this issue by noting that, for both Ribaucour transforms, the images
of the eigendirections of $\r_0$ and $\r_1$ under the
parametrisations $\phi_0,\tilde\phi_0$ of $s$ are tangent to systems
of circles on $s$ passing through $p_0$, $p_\infty$ or $\tilde
p_0$,$\tilde p_\infty$, respectively, together with their orthogonal
circles.  The two endomorphisms can only have the same
eigendirections (and so commute) if these two circle systems
coincide, that is, if
\begin{equation*}
\{\tilde p_0,\tilde p_\infty\}
   = \{p_0,p_\infty\}.
\end{equation*}
Thus, in the generic case where this condition is not met, the \bpt\
fails.

On the other hand, if we choose $\tilde p_0=p_0$ and $\tilde
p_\infty=p_\infty$ then $\tilde\phi_0=\phi_0$, $[\r_0,\r_1]=0$, and
the \bpt\ holds.  To exhibit the Demoulin families, it only remains
to determine four parallel sections of $V=\fh_0\oplus\fh_1$.  In the
case at hand, this can be done by inspection
\footnote{$\hat\sigma_0:=(1-2u+\frac{u^2}{2})\hat\kappa_{01}+\half\sigma_0$
and $\hat\sigma_1:=(1-2u+\frac{u^2}{2})\hat\kappa_{12}+\half\sigma_1$
complement the sections $\sigma_0$ and $\sigma_1$ to give sections of
$V$ which are parallel up to a common scaling by
$1/(1-2u+\frac{u^2}{2})$, as one easily verifies.} and it is then a
matter of linear algebra (which can be delegated to a computer
algebra engine) to determine the Demoulin families.

It turns out that both families consist of Dupin cyclides apart from
two spheres $s$ and $\tilde s=\Span{t_1-3e-2(p_0+p_\infty)}$ in the
$\alpha$-family (containing $f_0$) and one sphere (also $\tilde s$
but parametrised differently) in the
$\beta$-family (containing $\fh_0$ and $\fh_1$).  After suitable
stereographic projection $s$ becomes a plane, say $z=0$, $\fh_0$
becomes a vertical line, say $x=y=0$, and all Dupin cyclides become
surfaces of revolution with that line as axis\footnote{In
particular, $\fh_1$ is a circle parallel to $f_0$.}: their meridian
curves in the $y=0$-plane are shown in Figure~\ref{fig:3}.

\providecommand{\bysame}{\leavevmode\hbox to3em{\hrulefill}\thinspace}
\providecommand{\MR}{\relax\ifhmode\unskip\space\fi MR }
\providecommand{\MRhref}[2]{%
  \href{http://www.ams.org/mathscinet-getitem?mr=#1}{#2}
}
\providecommand{\href}[2]{#2}

\end{document}